\documentclass{amsart}
\usepackage{amssymb}
\usepackage{amsthm}
\usepackage{graphics}
\newtheorem*{reftheorem}{Theorem}
\newtheorem{introprop}{Proposition}
\newtheorem{introtheorem}{Theorem}
\newtheorem{theorem}{Theorem}

\newtheorem*{cor}{Corollary}
\newtheorem{prop}{Proposition}[section]

\newtheorem*{filllemma}{Filling Lemma}

\newtheorem{example}[prop]{Example}

\newtheorem*{dl}{Deformation Lemma}
\newtheorem*{lal}{Linear algebra lemma}
\newtheorem*{gromest}{Gromov's Hopf invariant estimate}
\newtheorem*{simplex}{Simplex straightening}
\numberwithin{equation}{section}
\title{The Hopf invariant and simplex straightening}
\author{Larry Guth}
\address{Department of Mathematics, University of Toronto, 40 St. George St., Toronto, ON M5S 2E4, Canada}
\email{lguth@math.toronto.edu}

\begin{document}
\begin{abstract}
Let $M$ be a closed oriented 3-manifold which can be triangulated
with $N$ simplices.  We prove that any map from $M$ to a genus 2
surface has Hopf invariant at most $C^N$.  Let $X$ be
a closed oriented hyperbolic 3-manifold with injectivity radius
less than $\epsilon$ at one point.  If there is a degree non-zero
map from $M$ to $X$, then we prove that $\epsilon$ is at least
$C^{-N}$.

\end{abstract}

\maketitle

In the 1970's, Thurston invented simplex straightening.  Milnor and Thurston
used simplex straightening to bound the degrees of maps to hyperbolic
manifolds.  In this paper, we extend their method to estimate the Hopf
invariant instead of the degree.

Here is a degree estimate that Milnor and Thurston proved.

\begin{reftheorem} (Milnor, Thurston) Let $M$ be a closed oriented
n-manifold that can be triangulated by N simplices, and let X be
a closed oriented hyperbolic n-manifold with volume V.  If f is
any continuous map from M to X, the norm of the degree of f is
bounded by $C N / V$.
\end{reftheorem}

We prove a similar estimate for the Hopf invariant.

\begin{introtheorem} Let $M$ be a closed oriented 3-manifold which can be
triangulated with $N$ simplices.  Let $f$ be a continuous map
from $M$ to a closed oriented surface of genus 2.  If the Hopf
invariant of $f$ is defined, then its norm is less than $C^N$.
\end{introtheorem}

(The Hopf invariant can be defined for a map from an 
oriented 3-manifold to an oriented surface provided that
the pullback of the fundamental cohomology class of the
surface is zero.  See Section 1 for more details.)

Our bound for the Hopf invariant grows exponentially in $N$.  
We will construct examples where the Hopf invariant is greater
than $(1 + c)^N$ for a universal constant $c > 0$.

Using estimates related to Theorem 1, we give a new proof of
a theorem of Soma:

\begin{reftheorem} (Soma) If $M$ is a closed oriented 3-manifold,
then there are only finitely many closed oriented hyperbolic
3-manifolds admitting maps of non-zero degree from M.
\end{reftheorem}

Soma's theorem follows from the following quantitative estimate:

\begin{introtheorem} Let $M$ be a closed oriented 3-manifold that
can be triangulated by $N$ simplices.  Let $X$ be a closed oriented
hyperbolic 3-manifold with injectivity radius $\epsilon$.  If
there is a map of non-zero degree from $M$ to $X$, then $\epsilon$ is
at least $C^{-N}$.
\end{introtheorem}

Roughly speaking, Theorem 2 says that a closed oriented hyperbolic 
manifold with small injectivity radius at one point must be topologically complicated.

Trying to understand Soma's theorem was the main motivation for
the work in this paper.  Here is some context for Soma's
theorem.  In dimension greater than 3, the analogue of Soma's
theorem is true, and the proof is much easier.  In dimension
greater than 3, there are only finitely many closed hyperbolic
manifolds with volume less than a given bound.  Applying the
theorem of Milnor and Thurston, it follows immediately that a
given n-manifold $M$ admits maps of non-zero degree to only
finitely many hyperbolic n-manifolds.  This proof does not work
in dimension 3, because there are infinitely many closed
hyperbolic 3-manifolds with volume less than a given bound.

In dimension greater than 3, the analogue of Soma's theorem also
holds for non-orientable n-manifolds.  The proof is essentially
the same as in the orientable case.  In \cite{BW}, Boileau and Wang proved
that the
analogue of Soma's theorem is false for non-orientable 3-manifolds.
Our proof of Theorem 2 requires $X$ to be orientable because we use
the Hopf invariant of maps from $X$ to the 2-sphere.  (If $X$ is non-orientable,
then the Hopf invariant is only defined modulo 2, and our arguments break down.)

The proof of Theorem 1 has two ingredients.  First, we use the simplex straightening method to
homotope an arbitrary continuous map to a map with a bounded Lipschitz constant.
Second, we use a method of Gromov (see \cite{G2}) to bound the Hopf invariant
of a map in terms of its Lipschitz constant.

The proof of Theorem 2 also begins with simplex straightening.  An arbitrary map
$f: M \rightarrow X$ can be straightened to a map with Lipschitz constant 1.  There are many obstructions to finding a degree non-zero
map with controlled Lipschitz constant.  The total volume gives the simplest obstruction,
and the volume estimate leads to Milnor-Thurston degree estimate.  But there are
other obstructions to finding a 1-Lipschitz map of non-zero degree.  In this paper,
we explore a different obstruction, related to a new Riemannian invariant called the Hopf size.

The Hopf size of a closed oriented 3-manifold $(M, g)$ measures the largest Hopf invariant of an $L$-Lipschitz map to the unit 2-sphere, with an appropriate normalization depending on $L$.

$$\textrm{Hopf Size } (M^3, g) := \sup_L \sup_{\textrm{Lip }(F) \le L} |\textrm{Hopf }(F)| L^{-4}.$$

We prove that the Hopf size of a closed oriented hyperbolic 3-manifold $X$ controls the degrees of maps to $X$.

\begin{introprop} Let $M$ be a closed oriented 3-manifold triangulated by $N$ simplices.  Let $X$ be a closed oriented hyperbolic 3-manifold.  Let $f: M \rightarrow X$ be a continuous map.

$$| \textrm{Degree } (f) | \le C^N / \textrm{Hopf Size }(X).$$

\end{introprop}

Proposition 1 is a variation of the Milnor-Thurston degree estimate, using the Hopf size of $X$ instead of the volume of $X$.  Proposition 1 gives new information not covered by the Milnor-Thurston estimate because a hyperbolic 3-manifold may have volume less than 100 and arbitrarily large Hopf size.

To prove Theorem 2, we perform a classification of closed oriented hyperbolic 3-manifolds with small injectivity radius.  Most of these hyperbolic manifolds have large Hopf size.  The rest of them are topologically complicated in a different way, such as having a torsion subgroup of $H_1$ with large order.

In Section 1, we review simplex straightening, the Hopf invariant, and Gromov's estimate.
In Section 2, we prove Theorem 1.  In Section 3, we introduce the Hopf size and prove Proposition 1.  In Section 4, we prove Theorem 2.  In Section 5 we give an example to show that the exponential growth
in Theorem 1 is real.  In Section 6 we give a brief exposition of the Boileau-Wang theorem
for maps to non-orientable 3-manifolds.

Notation: We use $|y|$ for the mass of a real chain $y$ in a Riemannian manifold.   We use the constant $C$ to denote a positive
constant that may change from line to line.

This paper is based on a section of my thesis, and I would like
to thank my advisor Tom Mrowka for his help and support.

\section{Background}

In this section we review simplex straightening, the Hopf invariant, and Gromov's
Hopf estimate.

First we describe simplex straightening.  By a hyperbolic k-simplex, we mean
a (finite) k-dimensional geodesic simplex in hyperbolic k-space.  By Lobachevsky's theorem, the volume of a hyperbolic k-simplex is bounded by $C(k)$ for each $k \ge 2$.  A hyperbolic simplicial
complex is a simplicial complex where each simplex is given the metric of
a hyperbolic simplex, and the metrics agree in the obvious way.

\begin{simplex} (Thurston) Suppose that $M$ is a simplicial complex and
$X$ is a complete hyperbolic manifold, and that $F_0$ is a continuous map
from $M$ to $X$.  Then we can choose hyperbolic metrics on each simplex
of $M$, and we can homotope $F_0$ to a new map $F$ in such a way that $F$ is a local isometry.  In particular, the Lipschitz constant of $F$ is 1.
\end{simplex}

Next we review the Hopf invariant.  Suppose that $F: M \rightarrow \Sigma$ is
a map from an oriented 3-manifold $M$ to an oriented surface $\Sigma$.  Let
$\omega$ denote the fundamental cohomology class in $H^2(\Sigma, \mathbb{R})$.
We will define the Hopf invariant of $F$ provided that $F^*(\omega) = 0$.

Let $q$ be a generic point in $\Sigma$.  Let $y = f^{-1}(q) \subset M$.  The fiber
$y$ is an oriented 1-manifold.  Since $F^*(\omega) = 0$, it follows by Poincare duality
that the homology class of $y$ vanishes in $H_1(M, \mathbb{R})$.  
Let $z$ be a real 2-chain with $\partial y = z$.
If we restrict the map $F$ to $z$, we get a map from $(z, \partial z)$ to $(\Sigma, q)$.
This map has a well-defined degree, which is equal to the Hopf invariant of $F$.  (The push-forward of $z$ is a real 2-cycle in $\Sigma$, and the Hopf invariant of $F$ is the homology class
of this cycle.)

By standard arguments from differential topology, the Hopf invariant does not
depend on the choice of $z$ or the choice of $q$, and it is a homotopy invariant
of the map $F$.  See \cite{Gu} for more background.

Now we turn to Gromov's geometric estimate for the Hopf invariant.  The question
is as follows.  Suppose that a map $F: (M^3, g) \rightarrow (\Sigma^2, h)$ has Lipschitz
constant $L$.  In terms of $L$ and geometric information about the domain and range,
how can we bound the Hopf invariant of $F$?

(To get a feeling for this problem, let's quickly consider the simpler case of 
bounding the degree of a map $F: (M^n, g) \rightarrow (N^n, h)$.  The
norm of the degree is at most $L^n Vol(M) / Vol(N)$.  For the Hopf invariant,
we cannot get a bound knowing only the Lipschitz constant $L$, the volume
of $M$, and the area of $\Sigma$.  We need some more refined geometric
information.)

Gromov's bound depends on understanding isoperimetric inequalities in $(M^3, g)$.
We define $Iso(M,g)$ as the smallest number so that every null-homologous
real 1-cycle $y$ in $(M,g)$ bounds a 2-chain of mass at most
$Iso(M,g) |y|$.  Gromov bounded the Hopf invariant of $F$ in terms
of the Lipschitz constant $L$, the volume of $(M,g)$, the area of $(\Sigma, h)$,
and this isoperimetric constant:

\begin{gromest} (\cite{G2}, see also \cite{Gu})  Suppose that $F: 
(M^3, g) \rightarrow (\Sigma^2, h)$ has Lipschitz constant $L$.
Then the Hopf
invariant of $F$ has norm at most $Iso(M,g) Vol(M,g) Area(\Sigma,h)^{-2} L^4$.
\end{gromest}

\begin{proof} (For a more detailed proof, see \cite{Gu}.)

By the coara formula, we can choose a generic $q \in \Sigma$ so
that the length of $y$ is at most $L^2 Vol(M) / Area(\Sigma)$.  By
assumption, we can choose $z$ with mass at most $Iso(M) Length(y) \le 
Iso(M) L^2 Vol(M) Area(\Sigma)^{-1}$.
The map from $z$ to $\Sigma$ has Lipschitz constant $L$.  So
its degree has norm at most $L^2 |z| Area(\Sigma)^{-1}
\le L^4 Iso(M) Vol(M) Area(\Sigma)^{-2}$. \end{proof}

For example, the Hopf invariant of an $L$-Lipschitz
map from the unit 3-sphere to the unit 2-sphere has norm at
most $C L^4$.  The exponent 4 here is sharp.  For $L > 2$, there
are maps from the unit 3-sphere to the unit 2-sphere with Lipschitz
constant $L$ and Hopf invariant at least $c L^4$.  For more
information, see \cite{Gu} or Chapter 7 of \cite{G}.

In order to apply Gromov's inequality, one has to estimate $Iso(M,g)$
for the Riemannian metric at hand.  The isoperimetric constant $Iso(M,g)$
is not as easy to compute/estimate as a volume or a Lipschitz constant.  
For example, if $(M,g)$ is a 3-dimensional ellipse, then $Iso(M,g)$
turns out to be roughly the second longest principal axis of
the ellipse, as proven in \cite{Gu}.  Even this special case takes
some work.

\section{Estimates for the Hopf invariant}

In this section, we prove Theorem 1.

\begin{theorem}  Let $M$ be a closed oriented 3-manifold which can be
triangulated with $N$ simplices.  Let $\Sigma$ be a closed surface of
genus 2 with fundamental cohomology class $\omega \in H^2(\Sigma, \mathbb{R})$.
Let $f$ be a continuous map $f: M \rightarrow \Sigma$ with $f^*(\omega) = 0$.
Then the Hopf invariant of $f$ has norm at most $C^N$.
\end{theorem}

\proof  Let $\Sigma$ be the surface of genus 2.  We equip it with
a hyperbolic metric of area $4 \pi$.  Next we apply simplex straightening
to the map $f$.  We get a new map $F$, homotopic to $f$, and
we give $M$ the structure of a hyperbolic simplicial complex.  We can
think of this structure as a Riemannian metric $g$ on $M$.  (The metric
is continuous but not smooth across the lower-dimensional faces of the triangulation.  This minor lack of regularity is not important.)  Simplex straightening tells us that
the map $F: (M^3, g) \rightarrow (\Sigma^2, hyp)$ has Lipschitz
constant 1.

We will estimate the Hopf invariant of $F$ using a variant of Gromov's method.
The technical difficulty that we face is that the space of hyperbolic simplices is
not compact.  Hence the space of possible metrics $g$ on $M$ is not compact.
We will prove a uniform estimate for the Hopf invariant of $F$.

The volume $Vol(M,g)$ is uniformly bounded by $C N$.  If we had
a uniform bound on $Iso(M,g)$, then we could apply Gromov's
Hopf invariant estimate and get a bound on the Hopf invariant of
$f$.  But I don't know any uniform bound for $Iso(M,g)$.  I suspect
there is no uniform bound for $Iso(M,g)$.  We will get around
this problem by using a small variation of Gromov's idea.  The
key estimate is the following isoperimetric inequality.

\begin{filllemma} Suppose that $(M^3, g)$ is a hyperbolic simplicial complex
with $N$ simplices.  Suppose that $y$ is an integral 1-cycle in $M^3$,
and that $y$ is null-homologous in $H_1(M, \mathbb{R})$.
Suppose that the mass of $y$ is at most $B$.  Also,
suppose that $y$ is transverse to the 2-skeleton of $M$, and
that $y$ intersects the 2-skeleton in at most $B$ points.

Then $y$ bounds a real 2-chain $z$ with mass at most $C^N B$.

\end{filllemma}

Using the Filling Lemma, we finish the proof of the theorem.

By Lobachevsky's theorem about the volumes of hyperbolic
simplices, we know that the volume of $(M^3, g)$ is at most
$C N$, and the area of the 2-skeleton is also at most $C N$.

By the coarea formula, we can choose a generic point $q \in \Sigma$
so that $y = F^{-1}(q)$ has length at most $C N$ and the number
of intersections of $y$ with the 2-skeleton is at most $C N$.  By Poincar\'e
duality $y$ is null-homologous.  By the Filling Lemma, we can choose a real 
2-chain $z$
with mass at most $N C^N$ with boundary $y$.  Since $F$ has
Lipschitz constant 1, the degree of the map from $z$ to $\Sigma$
has norm at most $C N C^N$.  Adjusting $C$, we see that the
Hopf invariant of $F$ is at most $C^N$. \endproof

Now we turn to the proof of the Filling Lemma.  The main ingredient is a variant
of the Federer-Fleming deformation theorem.

\begin{dl} Suppose that $(M^3, g)$ is a hyperbolic simplicial complex
with $N$ simplices.  Suppose that $y$ is an integral
1-cycle in $M^3$.  Suppose that the length of $y$ is at most $B$ and
that the number of intersections of $y$ with the 2-skeleton of $M$
is also at most $B$.

Then $y$ is homologous to a 1-cycle $y'$ in the 1-skeleton of $M$.  If
$e_i$ denote the (oriented) edges of $M$, then $y' = \sum_i c_i e_i$
with $c_i \in \mathbb{Z}$.  Moreover, $y'$ obeys the estimate

$$\sum_i |c_i| < C B $$

and $y'$ is homologous to $y$ by a 2-chain of area at most $C B$.

\end{dl}

\begin{proof}  We consider the intersection of $y$ with each
open 3-simplex in $M$.  Topologically, each component of
the intersection is either a segment with endpoints in the boundary
of a 3-simplex or else a closed curve in the simplex.

Step 1. (Straightening $y$.)  We construct a cycle $y_2$ simpex by simplex
in the following way.  For each topological segment in $y \cap \Delta^3$, we
put the geodesic segment in $\Delta^3$ with the same endpoints.

Now $y_2$ still has length at most $B$ and it still intersects the
2-skeleton of $M$ in at most $B$ points.  Moreover, the standard
isoperimetric inequality in hyperbolic space implies that $y$ is
homologous to $y_2$ by a 2-chain of area at most $C B$.

Step 2. (Pushing $y_2$ to the 2-skeleton.)  We construct a 1-cycle
$y_3$ from $y_2$ as follows.  For each simplex $\Delta^3$, and
each segment of $y_2$ in $\Delta^3$, we choose a path in
the boundary of $\Delta^3$ connecting the endpoints of the segment.  The path
in the boundary consists of two geodesic segments.

The cycle $y_3$ may be much longer than $y_2$.  It lies in the 2-skeleton of $M$, and it consists of
at most $C B$ geodesic segments.  Each of these segments has
one endpoint in the 1-skeleton of $M$.

The homology from $y_2$ to $y_3$ consists of at most $C B$
hyperbolic triangles, and so it has area at most $C B$.

Step 3. (Pushing $y_3$ to the 1-skeleton.)  The cycle $y_3$
is a union of geodesic segments.  Each segment has one
endpoint in the 1-skeleton of $M$.  Topologically, $y_3$
is a union of circles.  Each circle is a chain of geodesic
segments in the 2-skeleton of $M$.  Every other endpoint
in the chain of segments lies in the 1-skeleton of $M$.

If $p$ is an endpoint in the chain of segments lying in the interior
of a 2-face of $M$, the choose a vertex $v(p)$ in the boundary of
this 2-face.  Next, suppose that $\sigma \subset y_3$ is a geodesic
segment from $q$ to $p$, where $q$ lies in the 1-skeleton of $M$.
We replace $\sigma$ by a path from $q$ to $v(p)$, lying in the boundary of the
2-face containing $p$.  This path consists of at most two geodesic segments
each in the 1-skeleton of $M$.

We let $y_4$ be the sum of the new chains.  We see that
$y_4$ is a union of at most $C B$ geodesic segments,
where each segment lies in one edge of the 1-skeleton.
Hence $y_4 = \sum a_i e_i$ where $e_i$ are the edges
of the 1-skeleton and $\sum |a_i| \le C B$.

Next we bound the area of the homology from $y_3$ to $y_4$.
Let $\sigma'$ be the other segment touching $p$, and suppose that
$\sigma'$ goes from $p$ to $q'$.  Adding $\sigma$, $\sigma'$,
and their replacements, we get a 1-cycle in a single 2-simpex of $M$,
consisting of at most 6 geodesic segments.  It can be filled by
a chain consisting of at most 4 geodesic triangles.  By Lobachevsky's
theorem, the filling has area at most $C$.  Since the total number
of segments in $y_3$ was at most $C B$, the homology from $y_3$
to $y_4$ has mass at most $C B$.  \end{proof}

Now we fill this 1-cycle using linear algebra.

Let $e_i$ denote the edges of $M$ and $f_j$ denote the faces of $M$.  We assign
them orientations in an arbitrary way.  Then the simplicial boundary map
is given by a matrix $\beta_{ij}$, where $\partial f_j = \sum \beta_{ij} e_i$.
The entries of $\beta_{ij}$ are each $\pm 1$ or $0$.
The number of rows and the number of columns of the matrix are both
at most $C N$.  Because each 2-simplex has only three 1-simplices in its boundary,
each column of $\beta$ has at most three non-zero entries. 

\begin{lal} If $y' = \sum c_i e_i $ is a null-homologous real 1-cycle, then we can find a real 2-chain
$z = \sum d_j f_j$ with $\partial z = y'$, and $\sum_j |d_j| \le C^N \sum_i
|c_i|$. 

\end{lal}

\begin{proof} Recall that the Hilbert-Schmidt
norm of a matrix with entries $\beta_{ij}$ is defined to be $(\sum
|\beta_{ij}|^2)^{1/2}$.  The Hilbert-Schmidt norm of the boundary
matrix $\beta$ is bounded by $(C N)^{1/2}$.

Let r be the rank of the boundary matrix $\beta$.  We can choose
orthogonal matrices $O_1$ and $O_2$ so that $O_1 \beta O_2$ is zero
outside of the top left-hand $r \times r$ sub-matrix.  We define $\beta'$
to be this $r \times r$ sub-matrix.   Because the rank of $\beta$ is $r$, one of the $r \times r$
sub-determinants of $\beta$ is not zero.  Because the entries of $\beta$ are
all integers, this sub-determinant is an integer, and so its norm
is at least 1.  The norm of the determinant of $\beta'$ is at least as
large as the norm of any $r \times r$ sub-determinant of $\beta$. 
Therefore, the determinant of $\beta'$ has norm at least 1.  The
Hilbert-Schmidt norm of $\beta'$ is equal to that of $\beta$, which is less
than $(C N)^{1/2}$.  Any $(r-1) \times (r-1)$ sub-matrix of $\beta '$
has smaller Hilbert-Schmidt norm.  The determinant of an
arbitrary matrix $a_{ij}$ is bounded by $\Pi_i (\sum_j
|a_{ij}|^2)^{1/2}$.  It follows from this estimate and our bound
on the Hilbert-Schmidt norm that every $(r-1) \times (r-1)$
sub-matrix of $\beta'$ has determinant bounded by $C^N$.  Since the
determinant of $\beta'$ has norm  at least 1, the norm of every entry in the
inverse of $\beta'$ is bounded by $C^N$.  So the norm of the inverse of
$\beta'$ is bounded by $C N C^N$.

Therefore, the exact 1-chain $\sum c_i e_i$ is equal to the
boundary of a 2-chain $\sum d_j f_j$, with $\sum |d_j|$
bounded by $C N^2 C^N \sum |c_i|$.  After increasing the constant
$C$, we can say that $\sum |d_j|$ is bounded by $C^N \sum |c_i|$. 
\end{proof}

Assembling all of our steps, we have produced a real 2-chain
$z$ with $\partial z =  y $ and mass at most $C^N B$.  This
finishes the proof of the Filling Lemma and hence the proof of Theorem 1.

\section{The Hopf size}

In this section we define a notion of size based on the Hopf invariant.

Let $(M^3, g)$ be a closed oriented Riemannian 3-manifold.  We consider
maps $F$ from $(M^3, g)$ to the unit 2-sphere.  Let $\omega \in H^2(S^2, \mathbb{R})$
be the fundamental cohomology class of the 2-sphere.  We recall from Section 1 that
the Hopf invariant of $F$ is defined whenever $F^*(\omega) = 0$ in $H^2(M, \mathbb{R})$.
The Hopf size of $M$ is defined as follows.

$$\textrm{Hopf Size } (M^3, g) := \sup_L \sup_{\textrm{Lip }(F) \le L} |\textrm{Hopf }(F)| L^{-4}.$$

Here the second $\sup$ is taken over all maps $F$ from $(M^3, g)$ to the unit 2-sphere with Lipschitz constant at most $L$ and with $F^*(\omega) = 0$.  The Hopf size is a finite positive number.  One can
give an upper bound by using Gromov's Hopf invariant estimate from Section 1 - this shows that
the Hopf size is finite.  On the other hand, there are examples of maps from the unit 3-sphere to
the unit 2-sphere with Hopf invariant at least $c L^4$, which shows that the Hopf invariant of the
unit 3-sphere is positive.  But any closed oriented 3-manifold admits a Lipschitz degree 1 map to
the unit 3-sphere, and so it too has positive Hopf size, by Proposition 3.2 below.

The idea of Hopf size is essentially due to Gromov.  In \cite{G2}, Gromov gave a related definition for the `area' of a Riemannian 3-manifold.  Our definition is a small variation on his, adapted to the topology
problems in this paper.

In \cite{G}, Chapter 2, Gromov studied the analogous definition for the degree.  If $(M^n, g)$ is a closed
orientable manifold, then 

$$\textrm{Degree Size }(M^n, g) := \sup_L \sup_{\textrm{Lip }(F) \le L} |\textrm{Degree }(F)| L^{-n}.$$

Here the second $\sup$ is taken over all maps $F$ from $(M^n, g)$ to the unit n-sphere with Lipschitz
constant at most $L$.  Gromov proved that the degree size of $(M^n, g)$ agrees with the volume of $(M^n, g)$ up to a constant factor $C(n)$.  So the degree size essentially is the volume.  The Hopf size
is analogous to the degree size but subtler.  We will see below that it is not proportional to the volume.

The methods in the last section bound the Hopf size of a manifold triangulated into hyperbolic simplices.

\begin{prop} Let $(M^3, g)$ be a closed oriented manifold triangulated into $N$ simplices, so that the restriction of $g$ to each simplex is the metric of a hyperbolic simplex.

$$\textrm{Hopf Size} (M^3, g) \le C^N.$$

\end{prop}

\proof Consider a map $F$ from $(M^3, g)$ to the unit 2-sphere with Lipschitz constant $L$ and with
$F^*(\omega) = 0$.  By the coarea formula, we can choose $q \in S^2$ so that $F^{-1}(q)$ has length at most $C N L^2$ and meets
the 2-skeleton of $M^3$ in at most $C N L^2$ points.  By the Filling Lemma from Section 2, $F^{-1}(q)$ bounds a real 2-chain with mass at most $C^N C N L^2$.  But then the push-forward of $z$ to $S^2$ has mass at most
$C^N C N L^4$.  After redefining the constant $C$, the Hopf invariant of $F$ has norm at most $C^N L^4$.  \endproof

The Hopf size can control the degrees of Lipschitz maps between Riemannian manifolds.

\begin{prop} Suppose that $(M, g)$ and $(M', g')$ are closed oriented Riemannian 3-manifolds, and that $F$ is a map from $M$ to $M'$ with Lipschitz constant $L$.

$$ | \textrm{Degree } (f) | \le L^4 \frac{ \textrm{Hopf Size } (M, g)}{  \textrm{Hopf Size }(M' , g')} . $$

\end{prop}

\proof By definition, we can find a map $h$ from $(M', g')$ to the unit 2-sphere with $|Hopf (h) |  Lip(h)^{-4}$ as close as we like to the Hopf size of $(M', g')$.  Next we consider the composition $h \circ f$ from $(M, g)$ to the unit 2-sphere.  It has Lipschitz constant at most $L \cdot Lip(h)$ and Hopf invariant
equal to $deg(f) Hopf(h)$.  \endproof

Now using simplex straightening, the degree of a smooth map from $M$ to a hyperbolic 3-manifold $X$ can be bounded in terms of the Hopf size of $X$.

\begin{prop} Let $M$ be a closed oriented 3-manifold triangulated by $N$ simplices.  Let $X$ be a closed oriented hyperbolic 3-manifold.  Let $f: M \rightarrow X$ be a continuous map.

$$| \textrm{Degree } (f) | \le C^N / \textrm{Hopf Size }(X).$$

\end{prop}

\proof By simplex straightening, we can find a hyperbolic triangulation $g$ for $M$ and a map $F$ homotopic to $f$ with Lipschitz constant 1.  By Proposition 3.1., the Hopf size of $(M, g)$ is at most $C^N$.  By Proposition 3.2., the result follows. \endproof

\section{The proof of Theorem 2}

First we recall the statement of Theorem 2.

\begin{theorem} Let $M$ be a closed oriented 3-manifold that
can be triangulated by $N$ simplices.  Let $X$ be a closed oriented
hyperbolic 3-manifold with injectivity radius $\epsilon$.  If
there is a map of non-zero degree from $M$ to $X$, then $\epsilon$ is
at least $C^{-N}$.
\end{theorem}

The theorem basically says that a hyperbolic 3-manifold with
small injectivity radius is topologically complicated.  
We will show that a closed oriented
hyperbolic 3-manifold with small injectivity radius is
large / complicated in one of four ways.  First, it may have a
large volume.  Second, it may have a large Hopf size.  Third, there may
be a torsion element in $H_1(X, \mathbb{Z})$ with a large order. 
Fourth, it may require surfaces of large genus to span $H_2(X,
\mathbb{Z})$.  Any of these four features allow us to bound the
degree of a map from $M$ to $X$ in terms of the number of
simplices of $M$.

An important example is the case of hyperbolic Dehn fillings of a
finite-volume oriented hyperbolic manifold with a single cusp. 
For motivation, we explain what happens in this case.

Let $X_0$ be a finite-volume oriented hyperbolic manifold with a
single cusp.  Let $a$ and $b$ be a basis for the homology of the
boundary torus of $X_0$, which we are going to Dehn fill.  Choose
the basis so that the homology class $a$ bounds in $X_0$, and so
that the intersection number of $a$ and $b$ is 1.  Let $X(m,n)$
be the manifold formed by Dehn filling the boundary torus along
the curve homologous to $m a + n b$, for relatively prime numbers
$m$ and $n$.  According to Thurston's theory, the manifold
$X(m,n)$ admits a hyperbolic structure for all but finitely many
choices of $(m,n)$.  These hyperbolic manifolds have uniformly
bounded volume.  They contain a short core geodesic with length
on the order of $(m^2 + n^2)^{-1}$.  When $n=0$, this core
geodesic is non-torsion in $H_1(X(m,n), \mathbb{Z})$.  Otherwise,
the core geodesic is torsion with order $n$.  If $n$ is large,
then $X(m,n)$ is topologically complicated because of this
high-order torsion element.

Our new idea concerns what happens when $n$ is small but $m$
is large.  In this case, we will prove that there is a map from
$X(m,n)$ to the unit 2-sphere with Lipschitz constant $L$ and
with Hopf invariant of norm at least $\sim (m/n) L^4$.  Hence
$X(m,n)$ has Hopf size on the order of $|m/n|$.

The reason is that uniformly in $m$ and $n$, the manifold
$X(m,n)$ contains two disjoint thick tubes homologous to $b$,
which lie in the thick part of $X_0$.  Neither tube is
null-homologous in $X_0$, so their linking number is not defined
in $X_0$.  But in $X(m,n)$, the tubes are (rationally)
null-homologous, and a short calculation shows that their linking
number is equal to $\pm (m/n)$.  Linking numbers of tubes
are closely related to the Hopf invariant, and using these tubes,
we can construct a map to the unit 2-sphere with Lipschitz
constant bounded independent of $m,n$ and with Hopf invariant
at least $|m/n|$.

Now we turn to the formal argument in the general case.  Let 
$X$ be a closed oriented hyperbolic 3-manifold with injectivity
radius $\epsilon$.  Then X contains a closed embedded geodesic
$\gamma$ of length on the order of $\epsilon$.  If $\gamma$ is
torsion in homology with a fairly small order, then we will prove
that either X has large Hopf size or else X has large volume. 
This estimate is the main idea in the proof of Theorem 2.

\begin{prop} Let $X$ be a closed oriented hyperbolic 3-manifold
with a closed geodesic $\gamma$ of length $\epsilon$ which is
torsion in homology.  Then either the volume of $X$ is at least $c
\epsilon^{-1/6}$, or the order of $\gamma$ is at least $c
\epsilon^{-1/6}$, or the Hopf size of $X$ is at least $c
\epsilon^{-1/4}$.
\end{prop}

\proof Let T be the Margulis tube around $\gamma$.  (The thin
part of X is defined to be the subset of X where the injectivity
radius is less than a certain constant, and the Margulis
tube is the connected component of the thin region containing
$\gamma$.)  The universal cover of $\gamma$ is a geodesic in
hyperbolic 3-space.  We can parameterize hyperbolic space by the
upper half-space model so that the universal cover of $\gamma$ is
equal to the vertical line through the origin (i.e. the line
$x=0, y=0, z>0$).  The group of covering transformations of X
includes a loxodromic isometry that fixes this vertical line. 
This isometry is given by multiplying the three coordinates by a
constant on the order of $(1 + \epsilon)$, and rotating in the
(x,y)-plane by an angle $\theta$.  Consider the quotient of
hyperbolic 3-space by this isometry, and let $\tilde T$ be the
Margulis tube of this quotient around the core geodesic.  Taking
the quotient by the other covering transformations of X gives a
map from $\tilde T$ into M, which takes the core geodesic of
$\tilde T$ onto $\gamma$.  According to the Margulis Lemma, the
map is an embedding of $\tilde T$ into the Margulis tube around
$\gamma$.

The Margulis tube $\tilde T$ has a simple form.  For some number
R depending on $\epsilon$ and $\theta$, the tube $\tilde T$ is equal
to the quotient of the region $x^2 + y^2 < z^2 R^2$ by the action
of the loxodromic isometry corresponding to $\gamma$.  A
fundamental region for this action is given by the intersection
of the region above with the region $1 \le z \le 1+ \epsilon$. 
The boundary of the fundamental region includes two disks, the
first given by $z=1$ and $x^2 + y^2 \le R^2$, and the second
given by $z=1 + \epsilon$ and $x^2 + y^2 \le (1+\epsilon)^2 R^2$. 
The loxodromic isometry corresponding to $\gamma$ takes the first
disk onto the second disk.  In several constructions, we will use
the radial curves in $\tilde T$, which are the straight rays
through the origin in the Euclidean metric on the (x,y,z)-space. 
These radial lines are not geodesics in hyperbolic space, but
they are still useful in our proof.

We will give two estimates for R in terms of $\epsilon$ and
$\theta$.  The first estimate says that if $\epsilon$ is small,
then R is large.  More precisely, the radius R is at least $c
\epsilon^{-1/2}$.  Beginning at a point p on the edge of $\tilde
T$, we can follow a radial curve with length $N
\epsilon R$ going N times around the edge of $\tilde T$.  This
curve hits the circle $z=1$, $x^2 + y^2 = R^2$ in N points. 
Connecting the closest two of these N points by an arc of the
circle, we get a homotopically non-trivial closed curve with
length less than $N \epsilon R + R/N$.  Since this curve lies on
the edge of the Margulis tube, its length is greater than a
constant on the order 1.  We can make this construction for any
number N.  In particular, if $N = \epsilon^{-1/2}$, we get
the inequality $\epsilon^{1/2} R \ge c$.  This proves our lower
bound $R > c \epsilon^{-1/2}$.  On the other hand, we can assume
that R is not too big.  The volume of $\tilde T$ is roughly
$\epsilon R^2$.  If the volume of X is at least $c
\epsilon^{-1/6}$ then we are done, so we may assume that $R < C
\epsilon^{-7/12}$. 

Our second estimate says that if $\theta/(2\pi)$ is
well-approximated by rational numbers, then R is large. Again, we
begin at a point at radius R from the the center of the
horosphere $z=1$, and we trace a radial curve that goes
vertically q times around the tube, and then connect the endpoint
of this curve to the starting point within the horosphere $z=1$. 
The total length of this curve is roughly $q \epsilon R + |q
\theta/(2 \pi) - p| R$, where p is the integer that makes this
expression smallest.  Since this curve is homotopically
non-trivial and lies on the edge of the Margulis tube, it must
have length at least on the order of 1.  Therefore, $|\theta/(2
\pi) - p/q| > c (1/q) (1/R) - \epsilon$.

We choose generators for the homology of the boundary of $\tilde
T$, given by l, the longitude, and m the meridian.  (The choice
of longitude is related to the choice of $\theta$ in the
following way.  Take a radial curve on the edge of our
fundamental domain, going once around the tube from $z=1$ to $z=1
+ \epsilon$.  This line connects to a point on the base
horosphere $z=1$.  From that point, follow an arc of the circle
with directed length $- \theta$ to the initial point of the
radial curve.  The resulting closed curve is homologous to the
longitude.  We orient $\gamma$ so that the z coordinate increases
as we go along $\gamma$.  We orient the longitude to be
homologous to $\gamma$ and we orient the meridian so that
it has linking number 1 with $\gamma$.)  Let $p m + q l$ be the
homology class of a primitive curve which bounds in $X - T$. 
(The numbers p and q are relatively prime integers with $q \ge
0$.  Since the proposition assumes that $\gamma$ is torsion in
homology, q is not zero.  The number q is the order of $\gamma$
in $H_1(X)$.)  If q is at least $c \epsilon^{-1/6}$ we are done,
so we will assume that q is smaller than $c \epsilon^{-1/6}$. 
The inequality at the end of the last paragraph tells us that
$|\theta/(2 \pi) - p/q| > c
\epsilon^{3/4} - \epsilon$.  Our proposition is only interesting
for small values of $\epsilon$, so we may assume that $c
\epsilon^{3/4}$ dominates $\epsilon$, giving us the inequality
$|\theta/(2 \pi) - p/q| > c \epsilon^{3/4}$.

Finally, we estimate the Hopf size of X.  We claim that the
Hopf size of X is at least $c |\theta/(2 \pi) - p/q| R^2$,
which is greater than $c \epsilon^{-1/4}$.  To see this, we build
a collection of thin tubes in $\tilde T - \gamma$.  The tubes are
easiest to construct in case $\theta/(2 \pi)$ happens to be a
rational number $M/N$.  For the time being we make this
assumption.  Take a radial line in $\tilde T$ and follow it N
times around the tube $\tilde T$ until it comes back to its
initial point making a closed circle.  Then take a small
neighborhood of this circle, also consisting of a union of radial
lines.  We assume that the tube around our first radial line
meets the horosphere $z=1$ in a region of the following form,
using polar coordinates: $ |\theta - \theta_0| < \delta$, and $|r
- r_0| < r_0 \delta$, for some small number $\delta$.  (We assume
$r_0$ is bigger than 1, to avoid some unimportant technicalities.) 
The entire tube meets the horosphere $z=1$ in N disjoint regions
of this kind.  This tube is bilipshitz to the euclidean tube
$S^1(N r \epsilon) \times D^2(\delta)$ (with bilipshitz constant
20).  The tube meets the horosphere $z=1$ in a region of area at
most $N R \delta^2$.  This area bound shows that we can fit at
most $R \delta^{-2} N^{-1}$ disjoint tubes of this kind into
$\tilde T$.  The tubes fit together nicely, so this bound is
sharp.  To do the construction exactly, pick $\delta = \pi/N$. 
Since each tube has cross-sectional area $\delta^2$, the total
cross-sectional area of all these tubes is roughly $R/ N$.

Now we make the additional assumption that N is a multiple of q. 
This assumption implies that each of our tubes in null-homologous
in X.  Some of the tubes stay in the central part of $\tilde T$
given by $x^2 + y^2 < (1/4) z^2 R^2$.  We call these inner tubes. 
Other tubes are disjoint from this central part, and we call them
outer tubes.  (A few tubes lie partly in the central part, and we
throw them out.)  We now compute the linking number in X of an
inner tube with an outer tube.  The core radial curve of the
outer tube can be homotoped to the boundary of the Margulis tube
$T$ without crossing the inner tube, and it is homologous to $M m
+ N l$.  The homology class $p m + q l$ bounds in the complement
of T by assumption.  Therefore, $(N/q) p m + N l$ bounds in X. 
Subtracting such a boundary from our curve leaves a 1-cycle in
the boundary of T homologous to $(M - N p/q) m$.  Since the core
circle of the inner tube winds N times around T in the positive
direction, the linking number of m and the inner tube is N. 
Therefore, the linking number of our curve with the inner tube is
$(M - N p/q)N$.  We write this expression as $(M/N - p/q)N^2$. 
Since we assumed $\theta/(2 \pi) = M/N$, we can rewrite this
expression as $(\theta/(2 \pi) - p/q) N^2$.

Using these linked tubes, we construct Lipschitz maps from $X$ to $S^2$
with large Hopf invariants.  Each of our tubes is bilipschitz to
the product $S^1(N r \epsilon) \times D^2(\delta)$.  We pick a compactly
supported degree 1 map $\phi_0$ from $D^2(\delta)$ to the unit 2-sphere,
with Lipschitz constant at most $10 \delta^{-1}$.  Then we extend this
map to our tube by first projecting the tube to the disk and then using $\phi_0$.

Now we consider three different maps from $X$ to the unit 2-sphere, all with
Lipschitz constant $C \delta^{-1}$.  For $F_1$ use the tube map for 
each inner tube, and send all other points to the basepoint of $S^2$.  For $F_2$
use the tube map for each outer tube, and send all other points to the base point
of $S^2$.  And for $F_3$, use the tube map for each inner and outer tube, and
send all other points to the basepoint of $S^2$.  Now the Hopf invariant
is defined for all three maps $F_1$, $F_2$, and $F_3$.  In terms of our linking data,
we can directly compute $Hopf(F_3) - Hopf(F_1) - Hopf(F_2)$.  This difference is given
by twice the number of inner tubes times the number of outer tubes times the linking number
of an inner tube with an outer tube.  The number of inner tubes is at least $c R \delta^{-2}
N^{-1}$.  The number of outer tubes is lower bounded by the same expression.  And
the linking number is at least $| \theta/ ( 2 \pi) - p/q | N^2$.  So the difference in
Hopf invariants is at least $c R^2 | \theta / (2 \pi) - p/q| \delta^{-4} \ge c \epsilon^{-1/4} \delta^{-4}$.
Therefore, the Hopf size of $X$ is at least $c \epsilon^{-1/4}$.

So far, we have assumed that $\theta/(2 \pi) = M/N$, where M and
N are relatively prime integers and q divides N.  This assumption
is probably false, but it is easy to approximate $\theta/(2
\pi)$ as closely as we like by such a fraction $M/N$.  For a very
good approximation, the Margulis tube $\tilde T$ is bilipshitz to
the Margulis tube $T'$ corresponding to a loxodromic isometry along a
core geodesic of length $\epsilon$ with a twist of angle $(2 \pi)
M/ N$.  We construct all of the tubes and maps in $T'$ as
above and then pull them back to $\tilde T$ by our near-isometry. 
The resulting maps give the same estimate for the Hopf size.
\endproof

To complement this proposition, we prove a different kind of
inequality for hyperbolic 3-manifolds with short geodesics that
are non-torsion in homology.

\begin{prop} If a closed oriented hyperbolic 3-manifold X has a
closed geodesic which is non-trivial in $H_1(X, \mathbb{Q})$ of
length $\epsilon$, then any surface with non-zero intersection
number with $\gamma$ must have both area and genus at least $c
\epsilon^{-1/2}$.
\end{prop}

\proof We consider as above the Margulis tube around the short
geodesic $\gamma$, which has radius at least $R =
c \epsilon^{-1/2}$.  It contains thin tubes, each going around $N$
times, with total cross-sectional area at least $R/N$.  Any
surface with a non-zero intersection number with $\gamma$ meets
each of these tubes at least N times, and so it has total area at
least $N (R/N) = R > c \epsilon^{-1/2}$.  By the Thurston simplex
straightening argument, the genus of such a surface must be at
least $c \epsilon^{-1/2}$ as well. \endproof

We define the spanning genus of a 3-dimensional simplicial
complex M to be the smallest genus G so that $H_2(M, \mathbb{Q})$
is spanned by surfaces of genus at most G.  We can rephrase the
last proposition in this language: if X is a closed oriented
hyperbolic 3-manifold with a closed geodesic of length $\epsilon$
which is not torsion in $H_1(X, \mathbb{Z})$, then the spanning
genus of X is at least $c \epsilon^{-1/2}$.

In the last section, we proved that the Hopf size of a complex
with N simplices is less than $C^N$.  Similarly, we will prove
that the spanning genus of a complex with N simplices is less
than $C^N$.

\begin{prop} Let M be a simplicial 3-complex with N simplices. 
Then the spanning genus of M is less than $C^N$.
\end{prop}

\proof We consider the boundary map from the free vector space of
2-simplices of M to the free vector space of 1-simplices of M. 
The boundary map is given by a matrix $M$, with dimensions less than
$C N$, and with each entry equal to $\pm 1$ or $0$.  Moreover,
since the boundary of each 2-simplex is only three 1-simplices,
each column of the matrix has only three non-zero entries.  Let r
be the rank of the matrix $M$.  Pick r 2-simplices $\Delta_i$ in X,
so that $M(\Delta_i)$ gives a basis for the image of $M$.  After
renumbering the 2-simplices, we suppose that these are $\Delta_1$
through $\Delta_r$.

Now for every other 2-simplex $\Delta_j$, we can express
$M(\Delta_j)$ as a sum $a_{i,j} M(\Delta_i)$.  (In this formula,
the index j runs from $r+1$ to the last 2-simplex, and the index
i runs from 1 to r.)  Let I be the restriction of $M$ to the span
of the first r simplices.  The linear map I is an isomorphism
from $\mathbb{R}^r$ onto the image of $M$.  This isomorphism has
Hilbert-Schmidt norm less than $C N^{1/2}$, and it has
determinant at least 1.  Therefore, its inverse has norm at most
$C^N$.  For a fixed j, the sum $a_{i,j} \Delta_i$ is simply
$I^{-1}(M(\Delta_j))$.  Since $M(\Delta_j)$ has length $\sqrt 3$,
the vector $a_{i,j}$ has length less than $C^N$.  In particular,
each coefficient $a_{i,j}$ is less than $C^N$.

The coefficients $a_{i,j}$ don't have to be integers, but they
are rational numbers with controlled denominators.  The vectors
$M(\Delta_j)$ are integral vectors in the image of $M$.  The
vectors $M(\Delta_i)$ span a lattice L of integral vectors in the
image of $M$.  Let D be the index of the lattice L inside all of
the integral vectors in the image of $M$.  In particular, $D
M(\Delta_j)$ lies in the integral span of the $M(\Delta_i)$, and
so $D a_{i,j}$ is an integer for every i and j.  Since the
determinant of I is bounded by $C^N$, D is also bounded by $C^N$.

Now the integral 2-cycles $D \Delta_j - \sum_i D a_{i,j}
\Delta_i$ span the kernel of the matrix M.  Each of these
2-cycles has less than $N C^N C^N$ 2-simplices, and so it lies in
the span of some surfaces with genus less than $C^N$.  \endproof

Over the course of the last proof, we have also bounded the size
of the torsion subgroup of $H_1(M, \mathbb{Z})$.  This subgroup
is equal to the integral vectors in the image of $M$ modulo the
lattice L, and so D is the order of the torsion subgroup of
$H_1(M, \mathbb{Z})$.  As we showed above, this order is less
than $C^N$.  By an analogous argument, the torsion subgroup of
$H^2(M, \mathbb{Z})$ has order less than $C^N$.

We can now give the proof of Theorem 2.

\newtheorem*{theorem2}{Theorem 2}

\begin{theorem2} Let M be a closed oriented 3-dimensional
manifold which can be triangulated by N simplices.  Let X
be a closed oriented hyperbolic 3-manifold with injectivity
radius $\epsilon$. If there is a map of non-zero degree from M to
X, then $\epsilon$ is at least $C^{-N}$. 
\end{theorem2}

\proof Since X has injectivity radius $\epsilon$, it must contain a
closed geodesic of length less than $2 \epsilon$.  According to
Propositions 4.1 and 4.2, either X has volume at least
$c \epsilon^{-1/6}$, or $H_1(X, \mathbb{Z})$ has a torsion element
with order at least $c \epsilon^{-1/6}$, or X has Hopf volume at
least $c \epsilon^{-1/4}$, or X has spanning genus at least $c
\epsilon^{-1/2}$.

If X has volume at least $c \epsilon^{-1/6}$, then according to
Thurston's straightening theorem, N is greater than $c
\epsilon^{-1/6}$, and we are done.

If $H_1(X, \mathbb{Z})$ has a torsion element of order T, greater
than $c \epsilon^{-1/6}$, then we proceed as follows.  Let
$\alpha$ be this torsion class, and let $\alpha^*$ be the
Poincare dual class in $H^2(X, \mathbb{Z})$.  The element
$f^*(\alpha^*)$ in $H^2(M, \mathbb{Z})$ must also be torsion with
some order S dividing T.  We claim that T divides the product of
S with the degree of f.  Let $\alpha_T$ be the image of $\alpha$
in $H_1(X, \mathbb{Z}_T)$.  Let $\alpha_T^*$ be the Poincare dual
class in $H^2(X, \mathbb{Z}_T)$.  The pullback $S
f^*(\alpha_T^*)$ also vanishes.  Since $H_0(X, \mathbb{Z})$ is
free, $H^1(X, \mathbb{Z}_T)$ is equal to $Hom(H_1(X, \mathbb{Z}),
\mathbb{Z}_T)$.  In particular, we can choose a class $a$ in
$H^1(X, \mathbb{Z}_T)$ so that $a(\alpha)=a(\alpha_T)=1$ modulo
T.  Therefore, the cup product $a \cup \alpha_T^*$ is equal to
the fundamental cohomology class $O_T$ in $H^3(X, \mathbb{Z}_T)$. 
We let $[M]_T$ denote the image of the fundamental homology class
of M in $H_3(M, \mathbb{Z}_T)$.  The degree of f modulo T is
given by the pairing $f^*(O_T) ([M]_T)$.  Finally, S times the
degree of f modulo T is given by $f^*(S O_T) ([M]_T)$.  But $f^*(S
O_T)$ is equal to $f^*(a) \cup S f^*(\alpha_T^*)$, which vanishes. 
Since the volume of X is at least on the order of 1, the degree
of f is at most $C N$ by Thurston's simplex straightening
argument.  Also, the number S is bounded by the order of the
torsion subgroup of $H^2(M, \mathbb{Z})$, which is bounded by
$C^N$.  Therefore $C^N$ must be greater than $c \epsilon^{-1/6}$,
and we are done.

If the Hopf size of X is at least $c \epsilon^{-1/4}$, then Proposition
3.3 implies that $C^N \ge c \epsilon^{-1/4}$ and we are done.

If the spanning genus of X is at least $c \epsilon^{-1/2}$, we
proceed as follows.  Since the degree of f is non-zero, and since
X obeys Poincare duality, the map $f_*$ from $H_2(M, \mathbb{Q})$
to $H_2(X, \mathbb{Q})$ must be surjective.  Therefore, the
spanning genus of M must be at least $c \epsilon^{-1/2}$.  But
Proposition 4.3 tells us that the spanning genus of M is less
than $C^N$, and we are done. \endproof

The exponential bound in this theorem cannot be improved to a
sub-exponential bound.  To show this, we will construct closed
oriented hyperbolic manifolds X which can be triangulated by N
simplices and which contains closed geodesics of length less than
$c^{-N}$, for a constant $c>1$.  Begin with a non-compact finite
volume hyperbolic 3-manifold $X_0$ with a single cusp.  We view
$X_0$ as a manifold with boundary and we triangulate it.  Let T
be the restriction of our triangulation to the torus boundary. 
Next, we pick an Anosov diffeomorphism $\Psi$ of the torus and a
triangulation of the mapping cylinder of this diffeomorphism,
which restricts to the triangulation T on each boundary
component.  Finally, we pick a curve $c_0$ in the boundary torus
and a triangulation of a solid torus which restricts to T on the
boundary and which Dehn fills the curve $c_0$.  To form the
3-manifold X, we glue together $X_0$, N copies of the mapping
cylinder, and one copy of the solid torus.  (At each gluing, we
glue one torus with triangulation T to another torus with the
same triangulation by using the identity map.)  The resulting
3-manifold can be triangulated with less than $C N$ simplices. 
It is diffeomorphic to the Dehn filling of $X_0$ along the curve
$\Psi^{-N}(c_0)$.  If we fix any basis of $H_1(T^2)$, then for
most choices of $c_0$, the coefficients of the homology class
$\Psi^{-N}_*( [c_0])$ grow exponentially with N.  By Thurston's
theory, almost all of these Dehn fillings are closed hyperbolic
3-manifolds, with uniformly bounded volume.

If the short geodesic in the core of the Dehn filling has length
$\epsilon$, then the radius of the horosphere cross-section of
the Margulis tube must be roughly $\epsilon^{-1/2}$, and so there
is a homologically non-trivial simple curve in the boundary of
the Margulis tube of length less than $\epsilon^{-1/2}$ which
bounds a disk in the Margulis tube.  For large values of N, we
have proved above that this length must be at least $c^N$, for a
constant $c>1$.  Therefore, the length of the shortest geodesic
is less than $c^{-N}$.

As a corollary of Theorem 2, we can give a new proof of the
following result of Teruhiko Soma.  Trying to understand Soma's
result was the main motivation for the work in this paper.

\begin{cor} (Soma, \cite{S}) Given any closed oriented
3-manifold M, there are at most finitely many closed oriented
hyperbolic 3-manifolds X with maps of non-zero degree from M to
X.
\end{cor}

\proof Suppose that M can be triangulated with N simplices.  If
there is a map of non-zero degree from M to X, then Theorem 2
tells us that the injectivity radius of X is at least $C^{-N}$. 
The estimate of Milnor and Thurston tells us that the volume of X
is at most $CN$.  By Cheeger's finiteness theorem, there are
only finitely many hyperbolic manifolds obeying these bounds. \endproof

\section{Appendix 1: An example with large Hopf invariant}

This exponential upper bound seemed very high to me at first.  We
will now construct examples to show that it cannot be improved to
a sub-exponential bound.

\begin{example} There exists a triangulation of $S^3$ with N
simplices so that after giving each simplex a unit Euclidean metric,
the manifold has Hopf size greater than $exp(c N)$.
\end{example}

\proof We start with a torus $T^2$, equipped with a triangulation T
and with a choice of basis for $H_1(T^2)$, called a and b, with
intersection pairing of a and b equal to 1.  Next we take the
product $T^2 \times [0,1]$, and we Dehn fill both boundary
components to get $S^3$.  More precisely, we Dehn fill the
component $T^2 \times \{ 0 \}$ in such a way that the homology
class a bounds a disk in the new solid torus, and we Dehn fill
the component $T^2 \times \{ 1 \}$ in such a way that the
homology class b bounds a disk in the new solid torus.  We pick
triangulations of the two solid tori extending T.  (So far
we have only constructed the simplest Heegaard splitting of
$S^3$.)

Next we will construct a triangulation of the central cylinder
$T^2 \times [0,1]$ which restricts to the triangulation T on each
boundary component.  We pick an Anosov diffeomorphism $\Psi$ of
$T^2$.  For instance, $\Psi$ might act on homology by $\Psi_*(a)
= 2 a + b$ and $\Psi_*(b) = a + b$.  We pick a triangulation of
the mapping cylinder of $\Psi$ which restricts to the
triangulation T on each boundary component.  Our triangulation of
$T^2 \times [0,1]$ consists of N copies of this triangulated
mapping cylinder laid end to end, followed by N copies of the
mirror image of the mapping cylinder.  The whole triangulation
involves less than $C N$ simplices.

We will show that this triangulated 3-sphere contains two thick
tubes with linking number at least $exp(c N)$.  The first tube
$T_1$ is localized near the middle of the long cylinder.  Using
the local coordinates for $T^2$, it is homologous to a.  The
second tube $T_2$ is the core of the Dehn filling of $T^2 \times
\{ 0 \}$.  The linking number of $T_1$ with $T_2$ is equal to the
pairing of $\Psi_*^N(b)$ with $a$, which grows exponentially with
N.

Each tube $T_i$ is diffeomorphic to a product $D^2 \times S^1$.
We define a map $f_i$ from $T_i$ to $S^2$ by projecting to $D^2$
and then taking a degree 1 map from $D^2$ to $S^2$, mapping a
neighborhood of the boundary of $D^2$ to the basepoint of $S^2$. 
We fix these maps so that they don't depend on N.  In particular,
they have uniformly bounded 2-dilation.  Each map takes the
boundary of $T_i$ to the basepoint of $S^2$.  We extend the map
$f_i$ to all of $S^3$ by mapping the complement of $T_i$ to the
basepoint of $S^2$.  We construct a map $f_3$ from $S^3$ to
$S^2$, whose restriction to each tube $T_i$ is equal to $f_i$ and
which takes the rest of $S^3$ to the basepoint.

The maps $f_i$ have Lipschitz constant at most $C$.  The difference of Hopf
invariants $Hopf (f_3) - Hopf(f_1) - Hopf(f_2)$ is given by twice
the linking number of the two tubes, which is at least $(1 + c)^N$ for some
universal $c > 0$.
\endproof

The triangulated 3-manifold above has a large Hopf size, but it
is homeomorphic to $S^3$.  Therefore, any map from it to a
surface of genus 2 is contractible and has Hopf invariant zero. 
Using the triangulation above to guide a surgical operation, we
can construct manifolds $M^3$ which can be triangulated by N
simplices and which admit maps to a surface of genus 2 with Hopf
invariant at least $exp(cN)$.

\begin{example} We will construct closed oriented 3-manifolds M
which can be triangulated by N simplices, and which admit maps to
a surface of genus 2 with Hopf invariant at least $exp(c N)$.
\end{example}

\proof The manifold M is constructed by doing a kind of surgery
on the above triangulation near the two tubes $T_1$ and $T_2$. 
For each tube, we do the following procedure.  Identify T with
$D^2 \times S^1$, and refine the triangulation so that the
boundary of $D^2 \times S^1$ lies in the 2-skeleton.  Next, we
let $\Sigma'$ be a surface of genus two with one boundary
component.  We cut out $D^2 \times S^1$ from our manifold and
glue in $\Sigma' \times S^1$.  The gluing map is the
diffeomorphism from $\partial \Sigma' \times S^1$ to $\partial
D^2 \times S^1$ given by the product of a diffeomorphism from
$\partial \Sigma'$ to $\partial D^2$ and the identity map on
$S^1$.  Finally, we extend our triangulation to $\Sigma' \times
S^1$.  Since this operation occurs locally in a neighborhood of
T, it adds only C extra simplices independent of N.  Applying
this procedure to both $T_1$ and $T_2$, we get our 3-manifold M.

Now we construct some maps from M to a closed surface of genus 2. 
Inside of M, we have two copies of $\Sigma' \times S^1$, which we
added with our two surgeries.  For each copy, we can construct a
map in the following way.  We map $\Sigma' \times S^1$ to
$\Sigma'$ by projecting to the first factor, and then compose
with a degree 1 map from $\Sigma'$ to a surface of genus 2,
taking the boundary of $\Sigma'$ to the base point of the target
surface.  The resulting map takes the boundary of $\Sigma' \times
S^1$ to the base point of $\Sigma$, and so we can extend the map
to all of M, mapping the rest of M to the basepoint of $\Sigma$. 
In this way, we construct two maps, called $f_1$ and $f_2$.  We
can also construct a third map $f_3$ by applying the above
construction to both copies of $\Sigma' \times S^1$.  As in the
previous example, $\textrm{Hopf} (f_3) - \textrm{Hopf} (f_1) -
\textrm{Hopf} (f_2)$ is greater than $exp(c N)$. \endproof

\section{Appendix 2: The situation for non-orientable manifolds}

Our proof of Soma's theorem does not extend to non-orientable
manifolds, because intersection numbers, linking numbers, and the
Hopf invariant are only defined modulo 2 in non-orientable
manifolds.  As Boileau and Wang discovered, the analogue of Soma's theorem is
false for non-orientable manifolds.

\begin{reftheorem} (Boileau and Wang) There is a closed non-orientable 3-manifold M
which admits degree 1 maps to infinitely many different closed
oriented hyperbolic 3-manifolds X.  (The degree is only defined
modulo 2.)
\end{reftheorem}

The proof of this result is a small modification argument of an argument of Boileau and Wang in
\cite{BW}.  For reference, we include a proof here.

\proof The main idea of the proof is a clever choice of the
domain M which is due to Wang.  We begin with a finite volume
hyperbolic manifold $X_0$ with a single cusp, which we view as a
manifold with a torus boundary.  Let us pick a basis for
$H_1(\partial X_0)$, with elements a and b, so that the
intersection product of a and b is equal to 1.  We then perform a
filling of $X_0$ roughly analogous to a Dehn filling but using a
Mobius band in place of a disk.  Let us make a precise
construction.  Let B denote the Mobius band.  The boundary of $B
\times S^1$ is a torus.  We pick a diffeomorphism $\phi$ of the
boundary of $X_0$ with the boundary of $B \times S^1$, so that
$\phi(\partial B)$ is homologous to a and $\phi(S^1)$ is
homologous to b.  We define M to be the result of gluing $X_0$ to
$B \times S^1$ with this diffeomorphism.

We will construct degree 1 maps from M to half of all the Dehn
fillings of $X_0$.  All but finitely many of these Dehn fillings
admit hyperbolic metrics, and they include infinitely many
different manifolds.  For every pair $(m,n)$ of relatively prime
integers, let $X(m,n)$ denote the Dehn filling of $X_0$ along a
curve homologous to $ma + nb$.  The manifold $X(m,n)$ consists of
the union of $X_0$ and a solid torus, glued along their
boundaries.  The map that we will construct from M to $X(m,n)$
takes $X_0$ to $X_0$ identically.  The restriction of the
identity map to the boundary gives a map from the boundary of $B
\times S^1$ to the boundary of $D^2 \times S^1$. 

We now investigate in what cases we can extend this map to a map
from $B \times S^1$ to $D^2 \times S^1$.  Since maps to $D^2$
automatically extend, it suffices to extend the map to $S^1$ from
the boundary of $B \times S^1$ to the interior.

The homotopy classes of maps from the boundary of $B \times S^1$
to $S^1$ are simply $H^1(T^2)$.  The cohomology class of our map
is $ n (\partial B)^* -m (S^1)^*$, where $(\partial B)^*$ is the
cocycle that takes the value 1 on $\partial B$ and the value 0 on
$S^1$.  Now B is homotopic to $S^1$.  Let c denote a circle in B
which is a deformation retract of B.  Then $H^1(B \times S^1)$ is
generated by $c^*$ and $(S^1)^*$.  The inclusion map of $\partial
B \times S^1$ into $B \times S^1$ induces a map on cohomology
taking $c^*$ to $2 (\partial B)^*$ and $(S^1)^*$ to $(S^1)^*$.  A
given map from the boundary to $S^1$ extends to the interior if
and only if the corresponding class in $H^1(\partial B \times
S^1)$ lies in the image of $H^1(B \times S^1)$.  This condition
is met if n is even.  Therefore, for every even n, there is a
degree 1 map from the non-orientable 3-manifold M to the
Dehn filling $X(m,n)$.  \endproof

\end{document}